\documentclass[12pt]{amsart}
\usepackage{amsmath,amsfonts,amssymb,amscd,amsthm,amsbsy,epsf}
\newtheorem{prop}{Proposition}
\newtheorem{thm}[prop]{Theorem}

\def\be{\begin{equation}}
\def\ee{\end{equation}}
\def\Cech{{\v Cech}{} }

\def\zed{{\mathbb Z}}

\def\qed{\hbox{\hskip 6pt\vrule width6pt height7pt depth1pt \hskip1pt}}

\def\Z{{\mathbb Z}}

\def\R{{\mathbb R}}
\def\eps{\epsilon}

\def\tm{{\tilde \mu}}
\def\vecv{{\vec v}}

\begin{document}

\title{Topology of (some) tiling spaces without finite local complexity}
\author{Natalie Priebe Frank and Lorenzo Sadun}
\begin{abstract}
   A basic assumption of tiling theory is that adjacent tiles can meet
   in only a finite number of ways, up to rigid motions.  However,
   there are many interesting tiling spaces that do not have this
   property.  They have ``fault lines", along which tiles can slide
   past one another.  We investigate the topology of a certain class of
   tiling spaces of this type.  We show that they can be written as
   inverse limits of CW complexes, and their \Cech cohomology is
   related to properties of the fault lines.
\end{abstract}
\address{Natalie Priebe Frank\\Department of Mathematics\\Vassar
   College\\Box 248\\Poughkeepsie, NY 12604} \email{nafrank@vassar.edu}
\address{Lorenzo Sadun, Department of Mathematics, The University of
   Texas at Austin, Austin, TX 78712} \email{sadun@math.utexas.edu}
\maketitle


\markboth{Frank and Sadun}{Topology of tiling spaces without FLC}

\medskip

\section{Introduction}
\label{intro}

In discussing tilings, a standard assumption is that tiles can meet  
in only
a finite number of ways, up to rigid motion.  Equivalently, for any  
radius
$R$, there are only a finite number of patches of radius $R$ in the
tiling, up to rigid motion.  This condition is called {\em finite local
complexity}, or FLC.
\footnote{Frequently an even stronger condition is applied, namely
that tiles can only meet in a finite number of ways {\em up to  
translation},
a condition that excludes tilings like the pinwheel \cite{pinwheel},
in which tiles
appear in an infinite number of orientations.}

In this paper we consider tilings that do not meet the FLC  
condition.  We
show that spaces of such tilings can be given a natural topology in  
which
they are compact.  Many of the techniques used for FLC tilings, such as
inverse limit constructions and cohomology calculations, can be
modified to handle non-FLC tilings. In particular, we work out a number
of simple examples, and prove theorems about a broader class of  
examples.

\subsection{Past work}

The first substantial work on non-FLC tilings was done by Kenyon
\cite{Kenyon}, who considered a substitution that was
combinatorially equivalent to $a \to \begin{pmatrix} a& a& a \cr a & a
   & a \cr a & a & a \end{pmatrix}$, but in which one of the columns
was shifted by an irrational distance relative to the other two.  (The
tile itself was not a polygon.  Rather, it had straight edges on the
left and right, while the top and bottom had fractal pieces that
looked like the devil's staircase).  Upon further
substitution, the shift between the columns would grow to an
arbitrarily large multiple of the original irrational distance.  In
the limit, an infinite line of tile edges would appear, along which
tiles could face one another in an infinite variety of ways.  Such a
line is an example of a ``fault line", which we now define.

Given a finite set of {\em prototiles} (usually assumed to be
polygonal, or perhaps closed topological disks), a {\em tiling} is a
covering of $\R^d$ by rigid motions of copies of these prototiles that
are only allowed to intersect on their boundaries.  A {\em tiling
   space} $X$ is any translation-invariant set of such tilings that is
closed under the third topology of section \ref{threetops}.  Let $T\in
X$ be a tiling containing an infinite line (or ray) $\ell$ of tile  
edges, and
let $\vecv$ be a unit vector parallel to $\ell$.
We say that $\ell$ is a {\em fault line} if for every $\eps > 0$  
there is
a tiling $T' \in X$ such that:
\begin{enumerate}
\item On one half- (or quarter-) plane with boundary $\ell$, $T' = T 
$, and
\item there is a $t$ with $0 < |t| \le \eps$ such that on the other  
half-
   (or quarter-) plane with boundary $\ell$ we have $T = T' + t \vec 
{v}$.
\end{enumerate}
If for every sufficiently small $\eps > 0$ it is possible to choose  
$t=\eps$,
we call $\ell$ a {\em regular fault line}.

Fault lines play a central role in any discussion of non-FLC tilings.
This is because Kenyon also proves \cite{Kenyon2} that if a tiling
that is made from a finite prototile set has an infinite number of
inequivalent two-tile patches, then those patches occur along straight
edges, or they occur along an entire circle of tile boundaries.  The
former case leads to fault lines, while the latter can only occur for
special prototile sets (and never, for instance, in primitive
substitution tilings \cite{Frank}).

In 1998, Sadun \cite{Sadun} proposed some generalizations of the
pinwheel tiling.  One such example (Til(1/2)) uses polygonal tiles
(two similar right triangles) that appear in infinitely many
orientations, and has a fault line.
Interestingly, this example meets the conditions of Goodman-Strauss'
matching rules theorem \cite{Chaim}.  To wit: there exists a finite
collection of tiles and a finite set of local matching rules such that
these tiles tile the plane, but only in a manner that is locally
equivalent to the generalized pinwheel.  That is, a finite set of
local rules forces a global hierarchical structure that in turn forces
infinite local complexity.

\begin{figure}
\vbox{\epsfxsize=4truein\epsfbox{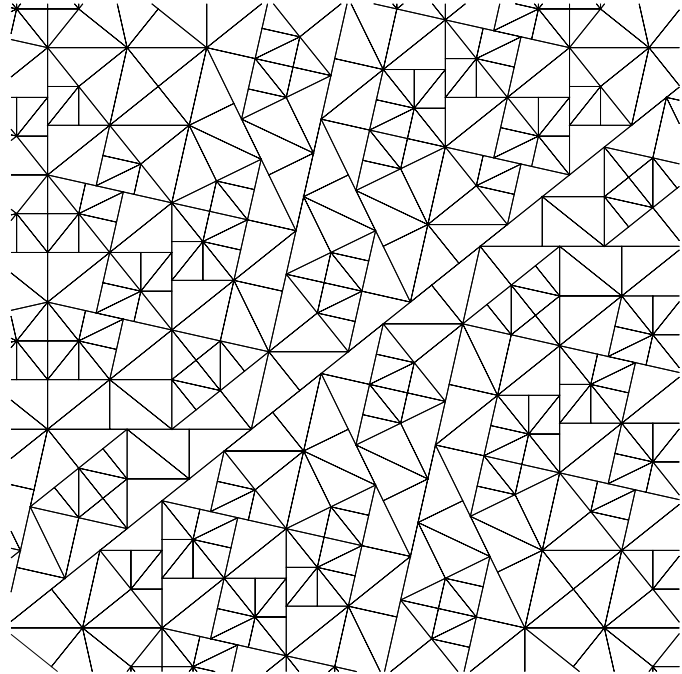}}
\caption{Part of the generalized pinwheel $Til(1/2)$.}
\end{figure}

Danzer \cite{Danzer} extended the theory of FLC tilings and provided
additional examples.  Finally, Frank and Robinson \cite{Frank2,  
Frank} have
considered a large family of ``direct product variation'' (DPV)
tilings.  These are obtained from products of 1-dimensional
substitutions by rearranging the positions of the tiles within an
order-1 supertile.  The examples of this paper are all DPV tilings in
which the rearrangements are purely horizontal, thereby preserving
the decomposition of the tiling into rows.

\begin{figure}
\vbox{\epsfxsize=4truein\epsfbox{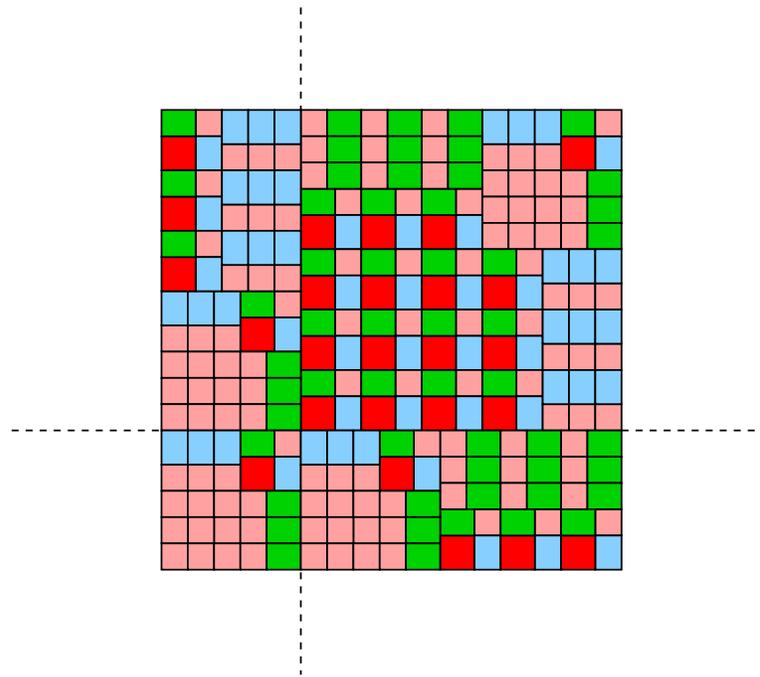}}
\caption{A patch of a DPV tiling with vertical and horizontal fault  
lines.}
\end{figure}

In all cases, the essential feature that prevents FLC is the presence of
a fault line.  Along the fault line, the evolution is described by two
1-dimensional substitutions.  One describes what happens along one  
side of
the fault line, while the other describes what happens along the  
other side.
The two substitutions cannot be Pisot with the same stretching  
factor, because in that case
the FLC property is preserved \cite{Frank} and there is no  
possibility of a fault line.
However, if the
substitutions are not Pisot, then differences in distributions of  
lengths
on opposite sides of the fault line will generally grow with successive
substitution, and the FLC property can be lost.  This will be discussed
further in section \ref{fault}.

\subsection{Three Topologies}
\label{threetops}

There are three metrics, and hence three topologies, that are frequently
applied to tilings and tiling spaces.

In the first metric, two tilings are considered $\epsilon$-close if they
agree, up to a translation of size $\epsilon$ or less, on a ball or
radius $1/\eps$ around the origin.  In this topology, the closure of the
(translational) orbit of a tiling is compact if and only if two  
conditions
are met: (a) there are only a finite number of tile types, up to  
translation,
and (b) tiles can only meet in a finite number of ways, up to  
translation.
This is the most frequently used topology for tilings that meet these
conditions.  Obviously, tiling spaces without FLC are not compact in
this topology, but neither are pinwheel-type spaces, which may have FLC
even though the tiles may appear in infinitely many orientations.

In the second metric, one applies a metric to the group of rigid motions
of the plane (e.g., defining an $\eps$ motion to be a translation of  
size
  $\le \eps$ followed by a rotation by $\le \eps$ about the origin). One
then considers two tilings to be $\eps$-close if they agree, up to an
$\eps$ motion, on a ball of size $1/\eps$ about the origin.  For tilings
in which tiles appear in only finitely many orientations (e.g., the  
Penrose
tiling), this is equivalent to the first topology.  However, it also  
handles
pinwheel-like spaces appropriately.  In this topology, the orbit closure
of a tiling is compact if and only if the tiling has FLC.

Finally, one can use a metric in which two tilings are $\eps$-close  
if they
contain the same tiles out to distance $1/\eps$, and if each tile in  
the first
tiling is within an $\eps$-motion of the corresponding tile in the
second tiling.  Note that in this topology, an $\eps$ shear along a  
fault
line yields a tiling that is $\eps$-close to the original, while in  
the first
two topologies it does not.

For tilings with FLC, the third topology is the same as the second,
insofar as it is impossible to apply a small rigid motion to one tile
without applying the same motion to all of its neighbors.  However, in
the third topology the orbit closure of any tiling with finitely many
tile types is always compact.  (To see sequential compactness,
start with an arbitrary sequence of tilings, pick a
subsequence in which the type, location and orientation of a tile near
the origin converges.  Then pick a subsequence in which the type,
location and orientation of a second tile converges.  Keep working
outwards from the origin, and then apply a Cantor diagonalization
argument to find a subsequence that converges everywhere.)

In the third topology, FLC is not a topologically invariant property.
Radin and Sadun \cite{RadSad} constructed a pair of spaces, one FLC
and one not, that are topologically conjugate.

\subsection{Outline of paper}

In section \ref{fault}, we study the evolution of fault lines in
2-dimensional substitution tilings.  This is essentially 1-dimensional
dynamics, and we relate properties of the fault line to the form
of the induced 1-dimensional substitutions.

In section \ref{simple}, we consider a 2-dimensional substitution tiling
with horizontal fault lines.  The rows of this tiling are (almost)
all the same, up to a horizontal shift, which is controlled by a  
vertical
1-dimensional substitution.  We show that the resulting tiling space can
be constructed as the inverse limit of compact CW complexes.  We  
explicitly
compute the cohomology of this tiling space, and discuss the meaning of
each term.

In section \ref{involved}, we consider a more complicated substitution,
as a step towards the direct product variations considered in section
\ref{general}.  These direct product variations look like the product of
a vertical and a horizontal 1-d substitution tiling, except that the  
rows
are sheared by an amount governed by the vertical substitution, and  
exhibit
horizontal fault lines.  The cohomology of the resulting tiling space is
computable in terms of the cohomologies of the vertical and  
horizontal 1-d
spaces, and the combinatorics of the vertical substitution.   
Specifically,
let $\mu$ be $H^1$ of the horizontal 1-d substitution tiling space,  
and let $M$
be the $n \times n$
substitution matrix of the vertical substitution, as applied to collared
tiles.  Then
\begin{itemize}
\begin{item} $H^0$ of the 2-d tiling space is $\zed$, of course. \end 
{item}
\begin{item} $H^1$ of the 2-d tiling space is isomorphic to $H^1$ of
the vertical 1-d substitution space. \end{item}
\begin{item} $H^2$ of the 2-d tiling space is isomorphic to the  
tensor product
of $\mu$ with the
direct limit of $\Z^n$ under the map $M^T$. This in turn is related to
$H^1$ of the vertical substitution space and the number of possible
fault lines. \end{item}
\begin{item} $H^3$ contains one copy of $\mu \otimes \mu$ for each  
possible
infinite fault line. \end{item}
\begin{item} $H^k$ is trivial for $k>3$. \end{item} \end{itemize}
Since $H^3$ is nontrivial, the 2-d tiling space is not homeomorphic  
to any
2-d tiling space with FLC.

Finally, in section \ref{open} we consider open problems in the  
theory of
tilings without FLC, and discuss our partial understanding of these  
problems.

\section{Analyzing a fault line --- 1 dimensional dynamics}
\label{fault}

Consider the 1-dimensional substitutions $\sigma_1(a)=ba$,
$\sigma_1(b)=aaa$, $\sigma_2(a)=ab$, $\sigma_2(b)=aaa$.  Both
substitutions have substitution matrix $\begin{pmatrix} 1 & 3 \cr 1 &
   0 \end{pmatrix}$, with Perron-Frobenius eigenvalue $\lambda = (1 +
\sqrt{13})/2 \approx 2.3028$, the larger root of the equation
$\lambda^2 - \lambda - 3 = 0$.  For a self-similar tiling, the $a$
tile can be given length $\lambda$ while the $b$ tile has length 3.
Note that for any word $W$, $\sigma_2(W)$ is a cyclic permutation of
$\sigma_1(W)$, obtained by removing an $a$ from the end and sticking
it on the beginning.  If $W$ is a bi-infinite word, then $\sigma_1(W)$
and $\sigma_2(W)$ are the same, up to translation by the length of
$a$.  This implies that the tiling spaces defined by $\sigma_1$ and
$\sigma_2$ are exactly the same.

Suppose we have a tiling with rectangular tiles $a$ and $b$ of widths
$\lambda$ and $3$, respectively, and suppose that $\sigma_1$ acts as a
substitution on lower edges and $\sigma_2$ acts as a substitution on
upper edges.  Let us construct a horizontal fault line, where the
evolution above the line is governed by $\sigma_1$ and the evolution
below the line is governed by $\sigma_2$.  If at some stage there is a
pair of exactly aligned $a$ tiles, one above the line and one below,
then on successive substitutions we will see

\be
\begin{pmatrix} b a  \cr a b \end{pmatrix}, \quad
\begin{pmatrix}  a a a b  a \cr a b a a a \end{pmatrix}, \quad
\begin{pmatrix}b a b a b a a a a b a \cr  a b a a a a b a b a b
                  \end{pmatrix}, \quad
\begin{pmatrix}
              a a a b a a a a b a a a a b a b a b a b a a a a b a \cr
a b a a a a b a b a b a b a a a a b a a a a b a a a
\end{pmatrix}.
\ee
Note that in the first and third substitution, the $a$ tiles are found
more on the right of the top row and on the left of the bottom row,  
while
in the second and fourth substitutions they are found more on the left
of the top row and the right of the bottom row.  The reason is that the
difference between the number of $a$ tiles up to a certain point  
grows like
the second eigenvalue of the substitution matrix, namely $1-\lambda  
\approx
-1.3$.  As we continue to iterate, this {\em discrepancy} grows without
bound.  (Strictly speaking, the discrepancy gets multiplied by $1- 
\lambda$
each time and then adjusted by $O(1)$ edge effects.  Once the  
discrepancy
grows beyond a certain point, the edge effects are dominated by the
multiplicative factor of $1-\lambda$ and we have exponential growth  
in the
discrepancy as a function of the number of substitutions. )

Pick a point along the fault line.  If there are $m$ more $a$ tiles  
in the
top row than the bottom up to that point, then the left edges of the
tiles on the top row
will be offset by $\lambda m \pmod{3}$ relative to the left edges of the
tiles on the bottom row.  By continuity, the discrepancy takes on all  
integer
values between 0 and $m$ as we move from the left edge of the pattern  
to the
point in question.  Since $m$ is unbounded and $\lambda$ is irrational,
this means that the possible  offsets of tiles in the top and bottom  
rows
takes on a dense set of values in the limit of infinite substitution.
In fact the left endpoints of upper $a$ tiles are dense in the lower  
$b$ tiles, because the only way for the
discrepancy to grow from $m$ to $m+1$ is for an additional $a$ tile  
to appear
along the top with a $b$ tile below it.    Thus every possible  
adjacency between
an upper $a$ an a lower $b$ can occur in the orbit closure; by  
primitivity of the
substitutions this implies that any adjacency is possible between any  
upper and lower tiles.
this, then
Thus we have not just a fault line, but a regular fault line.

Note how the form of the fault line depends on the second eigenvalue
of the substitution matrix.  If the substitution were Pisot, then the
discrepancy in the number of any species of tile would be bounded,
and the offsets between tiles would take on only a finite number of  
values.
As a result, the FLC condition would be preserved.  For a detailed
proof that Pisot substitutions do not lead to fault lines,
see \cite{Frank}.

Finally, note that there are only two possibilities involving  
substitutions
on two letters.  If the discrepancy in the number of $a$ tiles grows  
without
bound, then we have a regular fault line.  If the discrepancy is  
bounded,
then we preserve FLC.  It is not known whether irregular fault lines,
in which the set of possible offsets is infinite but not dense, are  
possible in substitution tilings.
In any case, they would require more than two letters.

\section{A simple 2-dimensional example}
\label{simple}

Consider a 2-dimensional tiling with two rectangular tiles.  Both the
$A$ and $B$ tiles have height 1, but the $A$ tile has width $\lambda=
(1+\sqrt{13})/2$ and the $B$ tile has width 3.  We consider the
self-affine substitution $\Sigma(A)= \begin{pmatrix} A&B\cr B&A
\end{pmatrix}$, $\Sigma(B) = \begin{pmatrix} A&A&A\cr A&A&A
\end{pmatrix}$.  For any $n \in \mathbb N$, an {\em $n$-supertile} is
a collection of tiles of the form $\Sigma^n(A)$ or $\Sigma^n(B)$.  A
tiling is {\em allowed} by the substitution if each of its finite
patches of tiles can be found in some $n$-supertile.  The
smallest closed set (under the third topology) containing all allowed
tilings is called the {\em substitution tiling space $X_\Sigma$}.  We
will see that measure-theoretically, almost all of the tilings in
$X_\Sigma$ are allowed by the substitution, but that the ones that are
not are the ones that make the topology different than in the FLC
case.

Note that whenever
two supertiles meet along a horizontal boundary, applying $\Sigma$
changes the bottom of the top supertile by $\sigma_1$ and the top of the
bottom supertile by $\sigma_2$.  By the results of section \ref 
{fault}, the
substitution tiling space $X_\Sigma$
defined by $\Sigma$ exhibits horizontal regular fault lines.

Not every row is subject to arbitrary shears.  The rows themselves are
labeled by points in the dyadic integers, describing their hierarchy
in the vertical substitution.  The label in the $n$th spot is 0 if the
row is in the lower $(n-1)$-supertile of its $n$-supertile and 1 if  
it is
in the upper $(n-1)$-supertile. If the labels of two adjacent rows  
differ
only in the first digit, then the dyadic label of the upper row begins
with a 0 and the label of the lower row begins with a 1.  One can see
that the sequence of tiles in the two rows are identical, with the
upper row offset horizontally by $\lambda$. If they differ only in the
first two digits, then the dyadic label of the upper row begins with
01 and the label of the lower begins with 10, and the upper row is the
same as the lower row, but offset by $\lambda^2-\lambda$.  If they
differ only in the first three digits, the upper and lower labels
begin with 001 and 110 resp., and then the rows are offset by
$\lambda^3-\lambda^2-\lambda$.  If they differ in the first $n$ digits
(and agree thereafter), then the upper and lower labels begin with
$0^{n-1} 1$ and $1^{n-1}0$ resp., and they are offset by $\lambda^n -
\lambda^{n-1}- \cdots - \lambda$. (In general, one sees the new offset
as $\lambda$ times the previous offset, minus $\lambda$.)  However, in
some tilings there exists a row with dyadic label $1111\ldots$ and an
adjacent row above it with label $0000\ldots$. These rows do not have
to have the same sequence of tiles, and their offset is arbitrary.

Put another way, all tilings in the tiling space contain
horizontal lines separating identical rows of tiles, offset by  
arbitrarily
large amounts.  However, in a small set of tilings (corresponding to a
single orbit in the dyadic solenoid that describes the vertical
hierarchical structure) there exists a fault line in which the tiling
above the fault is unrelated to the tiling below the fault.  These
special tilings have measure zero with respect to all
translation-invariant measures, and hence have no effect on
measure-theoretic properties of the tiling space.

However, they have a tremendous effect on the topology of the tiling  
space.
Indeed, this part of the tiling space is 3-dimensional! To describe such
a tiling we must specify the height of the infinite fault line (a  
point in
$\R^1$), the row immediately above the fault line (a point in a 1- 
dimensional
tiling space) and the row immediately below the fault line (a point in a
1-dimensional tiling space).

\subsection{The approximant $L$}

We will show how $X_\Sigma$ is the inverse limit of an approximant $L$
under a bonding map induced by $\Sigma$ (which we will again call
$\Sigma$).  The CW complex $L$ is actually 4-dimensional , but
$\Sigma: L \to L$ is not onto, and $\Sigma(L)$ is only a 3-dimensional
subset of $L$.  We will see that the inverse limit of $L$ under
$\Sigma$ exhibits the right combination of 2- and 3-dimensional
elements.

Let $\sigma=\sigma_1$, and let $K$ be an Anderson-Putnam complex of
the 1-dimensional substitution
$\sigma$, obtained by using collared tiles with a sufficiently
large radius $D$.  In order to ensure that $\sigma_2$ is a shift of
$\sigma_1$, we pick $D > (\lambda+3)/(\lambda-1)$ so that
$\lambda D > D + \lambda + 3$.  The 1-dimensional tiling
space $X_\sigma$ is the inverse limit of $K$ under a bonding map  
induced by
$\sigma$.  We then let
\be L = K \times K \times K \times [0,1] / \sim, \qquad (x,y,z,0)\sim 
(w,x,y,1).
\ee

\begin{figure}
\vbox{\epsfxsize=5truein\epsfbox{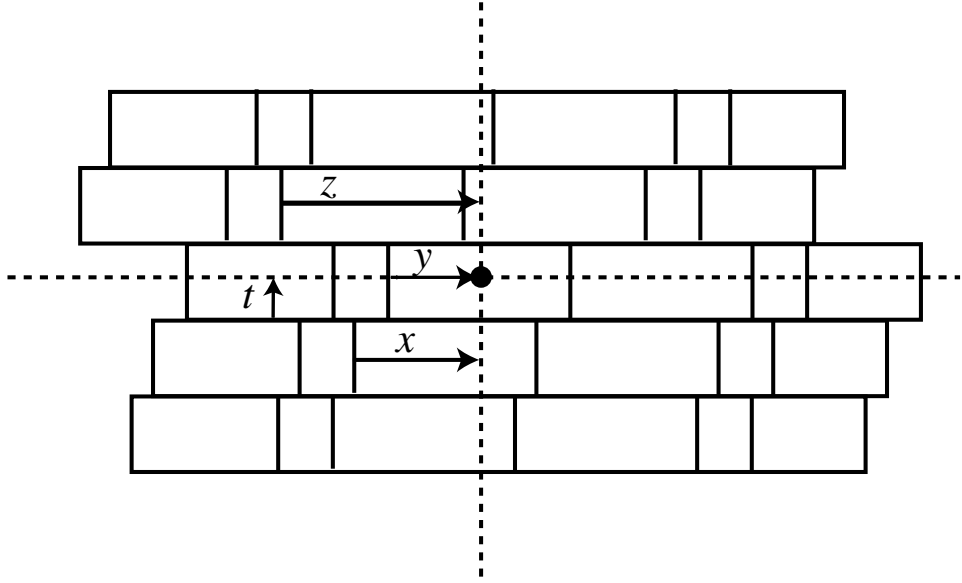}}
\caption{The data encoded in the approximant $L$.  The variables $x$,
$y$ and $z$ are points in $K$ that describe the rows nearest the origin,
while $t$ is the height of the origin in its row.}
\end{figure}

This is understood as follows.  Of the coordinates $(x,y,z,t)$ of a  
point
in $L$, $t$ describes the height of the origin in the row containing the
origin, and runs from 0 to 1.  The variable $y \in K$ describes a
horizontal neighborhood of size $D$ around the origin.  In other  
words, it
describes the row containing the origin.  The variables $x$ and $z$  
similarly
describe the rows immediately below and above, respectively.  If the  
origin
sits exactly on the boundary of two rows, we may describe the  
situation either
as $(x,y,z,0)$, with the origin sitting on the bottom of the ``$y$''  
row, or
as $(w,x,y,1)$, with the origin sitting on the top of the ``$x$''  
row. Under
the identification $\sim$, information about $w$ or $z$ is lost, and
we know only about the two rows touching the origin (``$x$'' and ``$y 
$'').

The bonding map induced from the substitution $\Sigma$ is
\be
\Sigma(x,y,z,t) = \begin{cases}
(\sigma_2(y), \sigma_1(y), \sigma_2(z), 2t-1); & \hbox{ if } t \ge  
1/2, \cr
(\sigma_1(x), \sigma_2(y), \sigma_1(y), 2t);  & \hbox{ if } t \le 1/2.
\end{cases}
\ee
However, $\sigma_2(y)$ is just a translate of $\sigma_1(y)=\sigma(y)$,
and likewise for $\sigma_2(z)$, so we can rewrite this as
\be \Sigma(x,y,z,t) = \begin{cases}
(\sigma(y)+\lambda,
   \sigma(y), \sigma(z)+\lambda, 2t-1); & \hbox{ if } t \ge 1/2, \cr
(\sigma(x), \sigma(y)+\lambda,
   \sigma(y), 2t); & \hbox{ if } t \le 1/2.
\end{cases}
\ee
Note that, depending on the value of $t$, the information from either  
$x$
or $z$ is lost, so that the image of $\Sigma$ is 3-dimensional.

Every translate is homotopic to the identity map, so the map $\Sigma$
is homotopic to
\be
\Sigma'(x,y,z,t) = \begin{cases}
(\sigma(y), \sigma(y), \sigma(z), 2t-1); & \hbox{ if } t \ge 1/2, \cr
(\sigma(x), \sigma(y), \sigma(y), 2t);  & \hbox{ if } t \le 1/2.
\end{cases}
\ee

\subsection{Computing $H^*(X_\Sigma)$}

To compute the cohomology of $X_\Sigma$, we simply compute the direct  
limit
of $H^*(L)$ under $\Sigma^* = (\Sigma')^*$.  Unlike $\Sigma$, $\Sigma'$
factors as a product of a vertical and horizontal map: $\Sigma'=
\Sigma_1 \circ \Sigma_2 = \Sigma_2 \circ \Sigma_1$, where
\be
\label{Sigma1}
\Sigma_1(x,y,z,t) = (\sigma(x), \sigma(y), \sigma(z), t), \ee
\be \label{Sigma2}
\Sigma_2(x,y,z,t) = \begin{cases}
(y,y,z,2t-1); & \hbox{ if } t \ge 1/2, \cr
(x,y,y,2t); & \hbox{ if } t \le 1/2.
\end{cases}
\ee

Since $\Sigma_1$ and $\Sigma_2$ commute, the direct limit of $H^*(L)$
under $\Sigma^* = \Sigma_1^* \circ \Sigma_2^*$ can be computed in two
steps.  First we take the direct limit under $\Sigma_1^*$, and then we
take the direct limit under $\Sigma_2^*$. Let $\tilde\mu=H^1(K)$, and
let $\mu=H^1(X_\sigma)$ be the direct limit of $\tilde\mu$ under $ 
\sigma^*$.
Our strategy is
\begin{enumerate}
\begin{item} Using Mayer-Vietoris, compute $H^*(L)$ in terms of $ 
\tilde\mu$.
\end{item}\begin{item}
Take the direct limit of $H^*(L)$ under $\Sigma_1$. Since $\Sigma_1$ is
essentially just the horizontal substitution $\sigma$, this merely
replaces each occurence of $\tilde\mu$ with $\mu$.  Note that we  
never have
to explicitly construct $K$ or compute $\tilde\mu$.
\end{item}\begin{item} Finally, take the direct limit under $ 
\Sigma_2^*$.
\end{item}\end{enumerate}

Step 1.  We take $V$ to be a neighborhood of $t=0$ (say, the set $t<0.2
\cup t>0.8$) and $U$ to be the region where $t$ is not close to zero
(say, $0.1<t<0.9$).  $U$ retracts to $K\times K \times K \times \{0.5 
\}$.
Let $\tilde \mu_x$ be the pullback of $\tilde \mu$ from the $x$  
factor, and
likewise for $\tilde \mu_y$ and $\tilde \mu_z$. We then have
\begin{eqnarray}
\notag H^0(U) & = & \Z, \\
\notag H^1(U) & = & \tilde \mu_x \oplus \tilde \mu_y \oplus
\tilde \mu_z, \\
H^2(U) & = & (\tm_x \otimes \tm_y) \oplus (\tm_x\otimes
\tm_z) \oplus (\tm_y \otimes \tm_z), \\
\notag H^3(U) & = & \tm_x \otimes \tm_y
\otimes \tm_z.
\end{eqnarray}
Likewise, $V$ retracts to $K \times K \times \{0\} \sim K \times K  
\times
\{1\}$. Let $\tm_{xy}$ and $\tm_{yz}$ denote the pullback of $\tm$  
from the
first  and second factors, respectively.
That is, $\tm_{xy}$ can be viewed either as $\tm_x$ from $t=0$ or
$\tm_y$ from $t=1$. We then have
\be H^0(V) = \Z, \qquad H^1(V) = \tm_{xy} \oplus \tm_{yz}, \qquad
H^2(V) = \tm_{xy} \otimes \tm_{yz}, \qquad H^3(V)=0. \ee
The intersection $U \cap V$ retracts to two copies of $K \times K  
\times K$,
say one at $t=0.15$ and at $t=0.85$, and we have
\begin{eqnarray}
\notag  H^0(U \cap V)& = & \Z^2, \\
\notag H^1(U\cap V) & = & (\tm_x \oplus \tm_y \oplus \tm_z) \oplus
(\tm_x \oplus \tm_y \oplus \tm_z), \\
H^2(U \cap V) & = &
((\tm_x \otimes \tm_y) \oplus (\tm_x\otimes \tm_z)
\oplus (\tm_y \otimes \tm_z)) \\
\notag && \oplus
((\tm_x \otimes \tm_y) \oplus (\tm_x\otimes \tm_z)
\oplus (\tm_y \otimes \tm_z)) \\
\notag H^3(U \cap V) & = &
(\tm_x \otimes \tm_y \otimes \tm_z) \oplus (\tm_x \otimes \tm_y
\otimes \tm_z).
\end{eqnarray}

The Mayer-Vietoris sequence is
\be
\cdots \longrightarrow H^k(L)
{\ \buildrel \rho \over \longrightarrow \ }
H^k(U)\oplus
H^k(V)
{\ \buildrel \nu \over \longrightarrow \ }
H^k(U \cap V)
{\ \buildrel \partial^* \over \longrightarrow \ }
H^{k+1}(L)
\longrightarrow
\cdots,
\ee
where $\rho$ is restriction and $\nu$ is signed restriction.  Using a  
basis for $\tm$ to make
bases for $H^k(U)$, $H^k(V)$ and $H^k(U\cap V)$,  and writing the  
lower copy of the basis of $H^k(U\cap V)$ first, the matrices of $\nu 
$ are
\be \nu = \begin{pmatrix} 1 & -1 \cr 1 & -1 \end{pmatrix} \qquad
\hbox{on $H^0$},
\ee
\be \nu = \begin{pmatrix} 1 & 0 & 0 & -1 & 0 \cr 0 & 1 & 0 & 0 & -1 \cr
0 & 0 & 1 & 0 & 0  \cr 1 & 0 & 0 & 0 & 0 \cr 0 & 1 & 0 & -1 & 0 \cr
0 & 0 & 1 & 0 & -1 \end{pmatrix} \qquad
\hbox{on $H^1$},
\ee
\be \nu = \begin{pmatrix}
1 & 0 & 0 & -1 \cr
0 & 1 & 0 & 0 \cr
0 & 0 & 1 & 0  \cr
1 & 0 & 0 & 0 \cr
0 & 1 & 0 & 0 & \cr
0 & 0 & 1 & -1
\end{pmatrix} \qquad
\hbox{on $H^2$},
\ee
\be \nu = \begin{pmatrix} 1 \cr 1 \end{pmatrix} \qquad
\hbox{on $H^3$},
\ee
where the $1$'s are actually identity matrices that depend on the  
dimension of $\tm$.
Note that these are all injective, except in dimension 0.  As a  
result, the
maps $\rho$ must all be zero (except in dimension 0), and for $k>0$ we
have that $H^k(L)$ is the cokernel of the previous $\nu$.  To summarize,
\begin{eqnarray}
\notag H^0(L) & = & \Z, \\
\notag H^1(L) & = & \Z, \\
H^2(L) & = & \tm, \\
\notag H^3(L) & = & (\tm\otimes\tm) \oplus (\tm\otimes \tm), \\
\notag H^4(L) & = & \tm\otimes\tm\otimes\tm.
\end{eqnarray}
By computing $\partial^*$ of the generators of $H^k(U \cap V)$, we  
see that
the generators of $H^*(L)$ are 1 in dimension 0,  $dt$ in dimension 1,
$dx\cup dt = dy \cup dt=dz\cup dt$ in dimension 2,
$dx\cup dy \cup dt = dy \cup dz \cup dt$ and $dx\cup dz\cup dt$ in
dimension 3 and $dx\cup dy \cup dz \cup dt$ in dimension 4. Here
we have used $dx$, $dy$, $dz$ and $dt$ as shorthand for 1-dimensional
cohomology generators in the $\tm_x$, $\tm_y$, $\tm_z$ and circle  
directions.

Step 2.  Taking the direct limit under $\Sigma_1^*$ merely converts
$\tm$ to $\mu$.  Using 1-dimensional methods \cite{BD} $\mu$
is easily shown to be the direct limit of $\Z^2$ under the matrix
$\begin{pmatrix} 1 & 1 \cr 3 & 0 \end{pmatrix}$, and is isomorphic to
$\Z[1/\lambda]$.

Step 3.  From equation (\ref{Sigma2}) we
compute the effect of $\Sigma_2^*$ on our cohomology
generators.
\begin{eqnarray}
\notag \Sigma_2^*(1) &=& 1 \\
\notag \Sigma_2^*(dt) &=& 2 dt \\
\Sigma_2^* (dy \cup dt) & =& 2 dy \cup dt \\
\notag \Sigma_2^* (dx \cup dy \cup dt) &=& dx \cup dy \cup dt \\
\notag \Sigma_2^* (dx \cup dz \cup dt) &=& dx \cup dy \cup dt
+ dy \cup dz\cup dt = 2 dx \cup dy \cup dt, \\
\notag \Sigma_2^*(dx \cup dy \cup dz \cup dt) &=& 0,
\end{eqnarray}
which gives us
\begin{eqnarray}
\notag H^0(X_\Sigma) & =& \Z, \\
\notag H^1(X_\Sigma) & = & \Z[1/2], \\
H^2(X_\Sigma) & = & \mu \otimes \Z[1/2] = \mu[1/2], \\
\notag H^3(X_\Sigma) & = & \mu \otimes \mu,
\\ \notag H^k(X_\Sigma) &=& 0 \hbox{ for } k>3.
\end{eqnarray}

Note that the identifications at the fault line prevent there being any
contribution of $\mu$ to $H^1(X_\Sigma)$.  We do get contributions from
$\mu$ to $H^2$ and $H^3$.  In particular, the $H^3$ term is easy to  
understand.
One factor of $\mu$ comes from the tiling above the fault line, one  
factor
of $\mu$ comes from the tiling below the fault line, and the $dt$ term
comes from the
location of the fault line. In this example, $H^2(X_\Sigma)$ equals the
tensor product of $\mu$ with $H^1$ of the dyadic solenoid; however,  
this is
not a general pattern. In the next section we shall see
an example in which $H^2$ of the 2-dimensional tiling space is not the
tensor product of $\mu$ with $H^1$ of the vertical substitution  
space, and in
section \ref{general} we will compute a general formula for $H^2(X_ 
\Sigma)$.

Finally, note that the form of the answer had nothing to do with the
details of the substitution $\sigma$, except that its expansion  
constant is not
a Pisot or Salem number.  Let $w_1$ and $w_2$ be any two words
in the letters $a$ and $b$, and consider the substitutions $\theta_1(a)=
w_1 a$, $\theta_2(b) =w_2a$, $\theta_2(a)=aw_1$, $\theta_2(b)=aw_2$. The
substitutions $\theta_1$ and $\theta_2$ generate the same 1-dimensional
tiling space.  If the pair $(\theta_1,\theta_2)$ generates a fault line,
as in section \ref{fault}, then consider the 2-dimensional substitution
$\Theta(A)=\begin{pmatrix} A & W_1 \cr W_1 & A \end{pmatrix}$,
$\Theta(B)=\begin{pmatrix} A & W_2 \cr W_2 & A \end{pmatrix}$, where
$W_1$ and $W_2$ are the same as $w_1$ and $w_2$, only written in capital
letters.  $\Theta$ gives rise
to a 2-dimensional tiling space with horizontal fault lines, and the
calculation of this section can be repeated, line by line, to show  
that the
cohomology of $X_{\Theta}$ is identical to that of $X_\Sigma$, only with
$\mu$ replaced by $H^1(X_{\theta})$.

\section{A more involved example}
\label{involved}

Our next example is a direct product variation, in the sense of
Frank \cite{Frank2}.
The vertical factor is the period-doubling substitution $0 \to 01$,  
$1 \to 00$,
while the horizontal factor is our usual 1-dimensional
substitution $\sigma$. We let
$A$ and $B$ denote $a \otimes 1$ and $b \otimes 1$, and abuse  
notation by
referring to $a \otimes 0$ and $b \otimes 0$ as $a$ and $b$,  
respectively.
The result is a tiling with four rectangular tiles. The tiles $A$ and  
$a$
have height 1 and width $\lambda$, while $B$ and $b$ have height 1 and
width 3.  Our substitution is
\be
\Sigma(a) = \begin{pmatrix} A&B \cr b&a \end{pmatrix}, \quad
\Sigma(A) = \begin{pmatrix} a&b \cr b&a \end{pmatrix}, \quad
\Sigma(b) = \begin{pmatrix} A&A&A \cr a&a&a \end{pmatrix}, \quad
\Sigma(B) = \begin{pmatrix} a&a&a \cr a&a&a \end{pmatrix}.
\ee
Note that if we ignore the difference between capital and lower case  
letters,
we revert to the example of section \ref{simple}.  The
period-doubling substitution space is an almost 1-1 extension of the  
dyadic
solenoid, and the example of this section is an almost 1-1 extension of
the example of section \ref{simple}.

\begin{figure}
\vbox{\epsfxsize=2.5truein\epsfbox{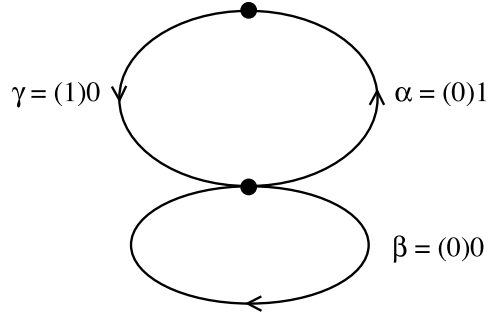}}
\caption{The Anderson-Putnam complex for the period-doubling  
substitution}
\label{AP-period}
\end{figure}

Before proceding with an analysis of $X_\Sigma$, we review some facts  
about
the period-doubling substitution.  This substitution forces the border
\cite{Kellendonk} on the right, since every substituted letter begins  
with a 0.
To construct the Anderson-Putnam complex \cite{AP}, we need only  
collar on the
left.  The complex (call it $P$, for period-doubling) is shown in  
figure \ref{AP-period}.
There are three collared tiles, which
we call $\alpha$, $\beta$ and $\gamma$.  The tile $\alpha$ is 1,  
preceded by
a 0, the tile $\beta$ is 0, preceded by a 0, and the tile $\gamma$ is 0,
preceded by a 1.  Viewed as a map on collared tiles, the substitution  
sends
$\alpha \to \gamma \beta$, $\beta \to \gamma \alpha$,
$\gamma \to \beta \alpha$, and interchanges the two vertices of the  
complex.

As before, let $K$ be the Anderson-Putnam complex for the horizontal
substitution space $X_\sigma$.  The approximant $L$ for $X_\Sigma$  
contains
a piece $K \times K \times K \times [0,1]$ for each edge $\alpha, \beta,
\gamma$ of $P$, with the identifications
\be (\alpha, x,y,z,0) \sim (\beta,x,y,z,0)\sim (\beta,w,x,y,1)\sim
(\gamma,w,x,y,1) \ee
and
\be (\gamma,x,y,z,0)\sim(\alpha,w,x,y,1).
\ee

The 5-tuple $(\alpha,x,y,z,t)$ (resp. $(\beta,x,y,z,t)$ or
$(\gamma,x,y,z,t)$) means that the origin is sitting at height $t$ in
a row that corresponds to $\alpha$ (resp. $\beta$ or $\gamma$) in the  
period-doubling
substitution.  Within the row containing the origin, the horizontal
position of the origin corresponds to $y \in K$.  In the rows below
and above that, the horizontal positions of the points one unit below  
and
above the origin correspond to $x \in K$ and $z \in K$.

The substitution $\Sigma$ induces the following map on $L$, also denoted
$\Sigma$:
\begin{eqnarray}
\notag \Sigma(\alpha,x,y,z,t) & = &
\begin{cases}
(\beta, \sigma_2(y), \sigma_1(y), \sigma_2(z), 2t-1); & \hbox{ if } t  
\ge 1/2,
\cr
(\gamma, \sigma_1(x), \sigma_2(y), \sigma_1(y), 2t);
& \hbox{ if } t \le 1/2.
\end{cases} \\
&=& \begin{cases}
(\beta, \sigma(y)+\lambda,
   \sigma(y), \sigma(z)+\lambda, 2t-1); & \hbox{ if } t \ge 1/2, \cr
(\gamma, \sigma(x), \sigma(y)+\lambda,
   \sigma(y), 2t); & \hbox{ if } t \le 1/2.  \end{cases} \\
\notag \Sigma(\beta,x,y,z,t) & = &
\begin{cases}
(\alpha, \sigma_2(y), \sigma_1(y), \sigma_2(z), 2t-1); & \hbox{ if }  
t \ge 1/2,
\cr
(\gamma, \sigma_1(x), \sigma_2(y), \sigma_1(y), 2t);
& \hbox{ if } t \le 1/2.
\end{cases} \\
&=& \begin{cases}
(\alpha, \sigma(y)+\lambda,
   \sigma(y), \sigma(z)+\lambda, 2t-1); & \hbox{ if } t \ge 1/2, \cr
(\gamma, \sigma(x), \sigma(y)+\lambda,
   \sigma(y), 2t); & \hbox{ if } t \le 1/2.
\end{cases} \\
\notag \Sigma(\gamma,x,y,z,t) & = &
\begin{cases}
(\alpha, \sigma_2(y), \sigma_1(y), \sigma_2(z), 2t-1); & \hbox{ if }  
t \ge 1/2,
\cr
(\beta, \sigma_1(x), \sigma_2(y), \sigma_1(y), 2t);
& \hbox{ if } t \le 1/2.
\end{cases} \\
&=& \begin{cases}
(\alpha, \sigma(y)+\lambda,
   \sigma(y), \sigma(z)+\lambda, 2t-1); & \hbox{ if } t \ge 1/2, \cr
(\beta, \sigma(x), \sigma(y)+\lambda,
   \sigma(y), 2t); & \hbox{ if } t \le 1/2.
\end{cases}
\end{eqnarray}

As before, $\Sigma$ is homotopic to a map $\Sigma'$ that does not  
contain
a translation by $\lambda$, and $\Sigma'=\Sigma_1\circ\Sigma_2 =
\Sigma_2 \circ \Sigma_1$, where $\Sigma_1$ implements $\sigma$  
horizontally,
and $\Sigma_2$ implements the period-doubling substitution vertically.
Specifically,
\be \Sigma_1(\zeta,x,y,z,t) = (\zeta,\sigma(x), \sigma(y), \sigma(z),  
t), \ee
where $\zeta = \alpha, \beta$ or $\gamma$, and
\begin{eqnarray}
\Sigma_2(\alpha,x,y,z,t) & = &
\begin{cases}
(\beta, y, y, z, 2t-1); & \hbox{ if } t \ge 1/2, \cr
(\gamma, x, y, y, 2t);
& \hbox{ if } t \le 1/2.
\end{cases} \\
\Sigma_2(\beta,x,y,z,t) & = &
\begin{cases}
(\alpha, y, y, z, 2t-1); & \hbox{ if } t \ge 1/2, \cr
(\gamma, x, y, y, 2t);
& \hbox{ if } t \le 1/2.
\end{cases} \\
\Sigma_2(\gamma,x,y,z,t) & = &
\begin{cases}
(\alpha, y, y, z, 2t-1); & \hbox{ if } t \ge 1/2, \cr
(\beta, x, y, y, 2t);
& \hbox{ if } t \le 1/2.
\end{cases}
\end{eqnarray}

We adopt the same strategy as in section \ref{simple}.  First we compute
the cohomology of $L$ using Mayer-Vietoris, then we take the direct  
limit
under $\Sigma_1^*$, and then we take the direct limit under $ 
\Sigma_2^*$.

To compute $H^*(L)$, let $V$ consist of neighborhoods of the two
vertices in $P$, and let $U$ be a slightly thickened complement to
$V$.  Now $U$ retracts to three copies of $K\times K\times K$, one for
each edge of $P$, and $V$ retracts to two copies of $K \times K$, one
for each vertex of $P$.  There are 6 copies of $K\times K \times K$ in
the retraction of $U \cap V$, one at the beginning of each edge and
one at the end of each edge.  We then have:
\begin{eqnarray}
   \notag H^0(U) & = & \Z_\alpha \oplus \Z_\beta \oplus \Z_\gamma, \\
   H^1(U) & = & \bigoplus_{e \in \{\alpha,\beta,\gamma\}}
(\tilde \mu_x \oplus \tilde \mu_y \oplus
   \tilde \mu_z)_e \\
   \notag H^2(U) & = & \bigoplus_{e \in \{\alpha,\beta,\gamma\}}
[(\tm_x \otimes \tm_y) \oplus (\tm_x\otimes
   \tm_z) \oplus (\tm_y \otimes \tm_z)]_e\\
\notag H^3(U) & = & \bigoplus_{e \in \{\alpha,\beta,\gamma\}}
(\tm_x \otimes \tm_y
\otimes \tm_z)_e.
\end{eqnarray}
\begin{eqnarray}
\notag H^0(V) & = & \Z^2, \\
H^1(V) & = & (\tm_{xy} \oplus \tm_{yz}) \oplus (\tm_{xy} \oplus \tm_ 
{yz})\\
\notag H^2(V) & = & (\tm_{xy} \otimes \tm_{yz}) \oplus
(\tm_{xy} \otimes \tm_{yz})\\
\notag H^3(V)& =& 0,
\end{eqnarray}
where the two copies in $H^k(V)$ correspond to the two vertices in $P$.
Since there is a copy of $H^k(U)$ at the beginning and at the end
of each edge of $H^k(U \cap V)$, we have that $H^k(U \cap V) = (H^k 
(U))^2$.

Writing down the generators of $H^k(U \cap V)$ in the same order
as the corresponding generators of $H^k(U)$ and $H^k(V)$ and keeping the
copies of $H^k(U)$ corresponding to the outgoing edges separate from  
those
corresponding to the incoming ones,
we can again compute matrices for $\nu$ in each dimension.   The  
matrices will
have two vertically aligned identity matrices in the $U$ columns and  
then
linearly independent vectors with two $-1$s in them in the $V$  
columns, so
the signed restriction maps $H^k(U) \oplus H^k(V) \to
H^k(U \cap V)$ are injective except in dimension zero. This implies that
the restriction maps $H^k(L) \to H^k(U) \oplus H^k(V)$ are all zero
(except in dimension zero), and that $H^{k+1}(L)$ is the cokernel of the
signed restriction map $H^k(U) \oplus H^k(V) \to H^k(U \cap V)$.
It is simple linear algebra to use the generators of $H^k(U \cap V)$ and
the image of $\nu$ to calculate that
\begin{eqnarray}
\notag H^0(L) & = & \Z, \\
\notag H^1(L) & = & \Z^2, \\
H^2(L) & = & \tm^5, \\
\notag H^3(L) & = & (\tm \otimes \tm)^7, \\
\notag H^4(L) & = & (\tm \otimes \tm \otimes \tm)^3,
\end{eqnarray}
with the following generators, subject to the following constraints.

Of course $H^0(L)$ is generated by 1.  The generators of $H^1(L)$ are  
$dt_\alpha$, $dt_\beta$
and $dt_\gamma$, with $dt_\alpha=dt_\gamma$. (In other words, $H^1(L) =
H^1(P)$.)  The products of generators of $\tm_x$, $\tm_y$ and
$\tm_z$ with $dt_\alpha$, $dt_\beta$ and $dt_\gamma$ generate $H^2(L) 
$.   Using the convention that
generators ranging through $\tm$ are referred to as simply $\tm$, we  
see the generators of $H^2(L)$ are
subject to the four constraints
\begin{eqnarray} \tm_y \cup dt_\alpha = \tm_x \cup dt_\gamma, &&
\tm_z \cup dt_\alpha = \tm_y \cup dt_\gamma, \\ \notag
\tm_y \cup (dt_\beta + dt_\gamma) = \tm_x \cup (dt_\alpha + dt_ 
\beta), &&
  \tm_z \cup (dt_\beta + dt_\gamma) = \tm_y \cup (dt_\alpha + dt_\beta).
\end{eqnarray}
The products of $(\tm_x \cup \tm_y)$, $(\tm_x
\cup \tm_z)$ and $(\tm_y \cup \tm_z)$ with $dt_\alpha$, $dt_\beta$
and $dt_\gamma$ generate
$H^3(L)$, subject to the two constraints
\be \tm_x\cup\tm_y\cup dt_\gamma = \tm_y \cup \tm_z \cup dt_\alpha,
\qquad \tm_x\cup \tm_y \cup(dt_\alpha + dt_\beta) = \tm_y \cup \tm_z  
\cup
(dt_\beta + dt_\gamma).
\ee
Finally, $H^4(L)$ is generated by $\tm_x\cup\tm_y\cup\tm_z\cup$ ($dt_ 
\alpha$,
$dt_\beta$, and $dt_\gamma$), with no constraints.

Note the form of the constraints on $H^2$ and $H^3$. For $H^2$, we  
have two
constraints for each vertex of $P$.  The sum of the $\tm_y \cup dt$  
terms
from the edges flowing into the vertex equals the sum of the $\tm_x  
\cup dt$
terms from the edges flowing out of the vertex, and the sum of the
$\tm_z \cup dt$ terms from the edges flowing into the vertex equals the
sum of the $\tm_y \cup dt$
terms from the edges flowing out of the vertex. These may be treated as
constraints among the $\tm_x \cup dt$ and $\tm_z \cup dt$ generators,  
while
the $\tm_y \cup dt$ generators are unconstrained.
For $H^3$ we have one
constraint per vertex, namely that the sum of the $\tm_y \cup \tm_z  
\cup dt$
terms from the edges flowing into the vertex equals the sum of the
$\tm_x \cup \tm_y \cup dt$ terms from the edges flowing out.

As in section \ref{simple}, the direct limit of $H^*(L)$ under $ 
\Sigma_1^*$
takes the same form as $H^*(L)$, only with $\tm$ replaced by $\mu$. What
remains is to take the direct limit under $\Sigma_2^*$, which we do one
dimension at a time.

The computation in dimension 0 is trivial, and we of course have
$H^0(X_\Sigma) =\Z$.

The computation in dimension 1 does not involve $\mu$ at all, and is  
identical
to the computation of $H^1$ of the period-doubling substitution space.
The answer is that $H^1(X_\Sigma)= \Z[1/2] \oplus \Z$.

For dimension 2, we first look at the $\mu_y \cup dt$ terms.  Since
$\Sigma_2^*$ maps these terms to themselves,
\begin{eqnarray}
\notag \Sigma_2^*(\mu_y \cup dt_\alpha) & = & \mu_y \cup dt_\beta +  
\mu_y \cup
dt_\gamma, \\
\Sigma_2^*(\mu_y \cup dt_\beta) & = & \mu_y \cup dt_\alpha + \mu_y \cup
dt_\gamma, \\
\notag \Sigma_2^*(\mu_y \cup dt_\gamma) & = & \mu_y \cup dt_\alpha +  
\mu_y \cup
dt_\beta,
\end{eqnarray}
these terms give the direct limit of $\Z^3$ under the matrix
$\begin{pmatrix} 0 & 1 & 1 \cr 1 & 0 & 1 \cr 1 & 1 & 0 \end{pmatrix}$.
This matrix has eigenvalues 2, -1 and -1, and the direct limit is $\Z 
[1/2]
\oplus \Z^2$.

Next we look at the $\mu_x \cup dt$ terms. We compute
\begin{eqnarray}
\notag \Sigma_2^*(\mu_x \cup dt_\alpha) & = & \mu_y \cup dt_\beta +  
\mu_y \cup
dt_\gamma, \\
\Sigma_2^*(\mu_x \cup dt_\beta) & = & \mu_y \cup dt_\alpha + \mu_x \cup
dt_\gamma = 2 \mu_y \cup dt_\beta, \\
\notag \Sigma_2^*(\mu_x \cup dt_\gamma) & = & \mu_x \cup dt_\alpha +  
\mu_x \cup
dt_\beta = \mu_y \cup dt_\beta + \mu_y \cup dt_\gamma.
\end{eqnarray}
Since the $\mu_x\cup dt$ terms map to the $\mu_y \cup dt$ terms, they do
not contribute anything additional to the direct limit.   Likewise, the
$\mu_z \cup dt$ terms map to the $\mu_y \cup dt$ terms, and do not  
give any
additional contributions.

In dimension 3, the image of $\Sigma_2^*$ consists entirely of
$\mu_x \cup \mu_y \cup dt$ and $\mu_y \cup \mu_z \cup dt$ terms,  
since the
image of $\Sigma_2$ never involves both $x$ and $z$.  Moreover, the $ 
\mu_x
\cup \mu_y \cup dt$ terms map to themselves, and the $\mu_y \cup
\mu_z \cup dt$ terms map to themselves.  Specifically, we have
\begin{eqnarray}
\notag \Sigma_2^*(\mu_x \cup \mu_y \cup dt_\alpha) & = & 0, \\
\notag \Sigma_2^*(\mu_x \cup \mu_y \cup dt_\beta)& = &\mu_x \cup \mu_y
\cup dt_\gamma, \\
\Sigma_2^*(\mu_x \cup \mu_y \cup dt_\gamma) & = & \mu_x \cup \mu_y \cup
(dt_\alpha + dt_\beta), \\
\notag \Sigma_2^*(\mu_y \cup \mu_z \cup dt_\alpha) & = & \mu_y \cup  
\mu_z
\cup (dt_\beta + dt_\gamma), \\
\notag \Sigma_2^*(\mu_y \cup \mu_z \cup dt_\beta)& = &\mu_y \cup \mu_z
\cup dt_\alpha, \\
\notag \Sigma_2^*(\mu_y \cup \mu_z \cup dt_\gamma) & = & 0.
\end{eqnarray}
Before taking the constraints into account, the action of $\Sigma_2^* 
$ has eigenvalues $1, 1, -1, -1, 0$ and $0$, making
the direct limit of this action $\Z^4$.  Because we know that
$\tm_x \cup \tm_y \cup dt_\gamma = \tm_y \cup \tm_z \cup dt_\alpha$  
and $\tm_x \cup \tm_y \cup
(dt_\alpha + dt_\beta) = \tm_y \cup \tm_z \cup (dt_\beta+ dt_\gamma)$  
in $H^3(L)$, and because $\Sigma_2^*$ annihilates
the $\mu_x \cup \mu_y \cup dt_\alpha$ and $\mu_y \cup \mu_z \cup dt_ 
\gamma$ terms, there are only two
independent generators left in the direct limit for $H^3(X_\Sigma) 
$.   These are $\mu_x \cup \mu_y \cup
dt_\gamma = \mu_y \cup \mu_z \cup t_\alpha$ and $\mu_x \cup \mu_y \cup
dt_\beta = \mu_y \cup \mu_z \cup dt_\beta$.    The cohomology has  
registered the two possible fault lines,
one between the $\alpha$ and $\gamma$ rows, and one between two $\beta 
$ rows.

In dimension 4, $\Sigma_2^*$ is identically zero.

To summarize, the cohomology of $X_\Sigma$ is
\begin{eqnarray}
\notag H^0 & = & \Z, \\
\notag H^1 & = & \Z[1/2] \oplus \Z = H^1(X_{\hbox{\small{pd}}}), \\
H^2 & = & \mu \oplus \mu \oplus \mu[1/2]
= \mu \otimes(H^1(X_{\hbox{\small{pd}}}) \oplus \Z), \\
\notag H^3 & = & (\mu \otimes \mu) \oplus (\mu \otimes \mu), \\
\notag H^k & = & 0 \hbox{ for }k>3.
\end{eqnarray}

Note that $H^2(X_\Sigma)$ is {\em not} the tensor product of $\mu$  
with the
first cohomology of the period-doubling substitution.  Rather, it
is the tensor product of $\mu$ with the direct limit of the transpose
of the period-doubling substitution matrix (as applied to collared  
tiles),
which has an additional factor of $\Z$.
The two copies of $\mu\otimes\mu$ in $H^3$ refer to the two types of
fault lines that can occur in $X_\Sigma$.  One has a $\beta$ row both  
above
and below the fault line.  The other has an $\alpha$ row below the  
fault line
and a $\gamma$ row above.  $\Sigma$ maps each of these situations to the
other.

\section{General direct product variations with regular fault lines.}
\label{general}

The results of the previous section are suggestive of how the fault  
lines
affect the cohomology of a tiling space.  In this section we prove  
two theorems
on the cohomology of tiling spaces with horizontal fault lines that  
arise as
direct product variations.

Suppose we have a collection of primitive
1-dimensional substitutions $\sigma_1$,
$\sigma_2, ..., \sigma_N$ defined on the same alphabet $A$.
Moreover, assume that for each $a \in A$, the length of $\sigma_k(a)$  
is the
same for all $k$, that all of the substitutions have the same  
stretching factor, and
that they all yield the
same tiling space.  (E.g., all of the $\sigma_k$s might be cyclic  
permutations
of one another).  Now suppose we have another primitive
1-dimensional substitution
$\rho$.  We could then consider the direct product of the two  
substitutions,
with $\sigma_1$ acting horizontally and $\rho$ acting vertically.  We  
then
replace $\sigma_1$ with $\sigma_2$, etc in various rows, so as to  
introduce
fault lines.  Furthermore, we assume that regular fault lines occur  
at every
boundary of infinite-order $\rho$-supertiles.

We call the resulting 2-dimensional substitution $\Sigma$, and  
compute the
cohomology of $X_\Sigma$.  Let $\mu = H^1(X_\sigma)$ and let $\nu=H^1 
(X_\rho)$.
Let $n$ be the number of configurations in the $X_\rho$ tiling space in
which two infinite-order $\rho$-supertiles meet at the origin.  By  
assumption,
this is the same as the number of ways that a regular fault line can  
occur in
an $X_\Sigma$ tiling.

\begin{thm}
Under these circumstances, the cohomology of $X_\Sigma$ is as follows:
\begin{eqnarray} \notag H^0(X_\Sigma)& = & \Z, \\
\notag H^1(X_\Sigma) & = & \nu, \\
H^2(X_\Sigma) & = & \mu \otimes (\nu \oplus \Z^{n-1}), \\
\notag H^3(X_\Sigma) & = & (\mu \otimes \mu)^n, \\
\notag H^k(X_\Sigma) & = & 0 \hbox{ for }k>3.
\end{eqnarray}
\end{thm}

\noindent Proof: As usual, let $P$ be the AP complex of
$X_\rho$, and let $K$ be the AP complex of $X_\sigma$.
We take $L$ to be one copy of $K \times K \times K \times [0,1]$ for  
each
edge in $P$, with the following identification.  If the end of edge $ 
\alpha$
meets the beginning of edge $\beta$ at a vertex in $P$, then
$(\alpha,w,x,y,1)\sim (\beta,x,y,z,0)$.

We compute $H^*(L)$ by Mayer-Vietoris.  Let $V$ be a union of  
neighborhoods
of the vertices of $P$, and let $U$ contain the edges.  $U$ retracts to
a number of copies of $K\times K \times K$ (one per edge),
while $V$ retracts to a number of copies of $K \times K$ (one per  
vertex).
$U \cap V$ retracts to a number of copies of $K \times K \times K$ (two
per edge, one at the beginning and one at the end).

Using the same bookkeeping as in section \ref{involved},  one can
see that the signed restriction maps $H^k(U) \otimes H^k(V) \to H^k(U 
\cap V)$ are
all injective, except in dimension 0.  The cokernel is understood as  
follows.
The image of the restriction of $H^k(U)$ to $H^k(U\cap V)$ merely
identifies (up to sign) the
two $K\times K \times K$ terms in $U \cap V$ that correspond to each  
edge.
This gives obvious generators for $H^k(L)$: In dimension 1 we have  
$dt_\zeta$,
where $\zeta$ is an edge in $P$, in dimension 2 we have $\tm_x \cup  
dt_\zeta$,
$\tm_y \cup dt_\zeta$, and $\tm_z \cup dt_\zeta$, in dimension 3 we have
$\tm_x \cup \tm_y \cup dt_\zeta$, $\tm_x \cup \tm_z \cup dt_\zeta$,
and $\tm_y \cup \tm_z \cup dt_\zeta$, and in dimension 4 we have
$\tm_x \cup \tm_y \cup \tm_z \cup dt_\zeta$, all subject to the  
identifications
imposed by the image of $H^k(V)$ in $H^k(U \cap V)$. In $H^1$ this is  
that
the sum of the $dt_\zeta$s entering a vertex equals the sum of the  
$dt_\zeta$s
coming out.  In $H^2$ it is that the sum of the $\tm_y \cup dt$ terms
from the edges flowing into the vertex equals the sum of the $\tm_x  
\cup dt$
terms from the edges flowing out of the vertex, and the sum of the
$\tm_z \cup dt$ terms from the edges flowing into the vertex equals the
sum of the $\tm_y \cup dt$
terms from the edges flowing out of the vertex.
For $H^3$, the the sum of the $\tm_y \cup \tm_z \cup dt$
terms from the edges flowing into a vertex equals the sum of the
$\tm_x \cup \tm_y \cup dt$ terms from the edges flowing out.
For $H^4$, there are no constraints.  This completes the computation of
$H^k(L)$.

Now we note that $\Sigma$ is homotopic to the product of $\sigma$ in the
horizontal direction and $\rho$ in the vertical direction.  Taking
the direct limit under $\sigma^*$ merely converts $\tm$ to $\mu$.  
What is
left is taking the direct limit under $\rho^*$.

In dimension 1, this gives $\nu$.

In dimension 2, we note that the constraints express certain  
combinations
of the $\mu_x \cup dt_\zeta$s or the $\mu_z\cup dt_\zeta$s
in terms of the $\mu_y \cup dt_\zeta$s.  They do not constrain the
$\mu_y \cup dt_\zeta$s terms, which are mapped to themselves
by $\rho^*$.  Furthermore, among these
terms the action of $\rho^*$
is just the transpose of the substitution matrix
itself (as applied to edges of $P$, i.e. to collared tiles).  The direct
limit of $H^2(L)$ under $\rho^*$ therefore contains the direct limit of
this matrix.

In principle, the direct limit of $H^2(L)$ should also
contain contributions from the $\mu_x\cup dt_\zeta$ and
$\mu_z \cup dt_\zeta$ terms.  However, we claim that
these contributions are zero as
a consequence of our using a complex that forces the border, as
explained below.

Note that the pullback of $dx \cup dt_\zeta$ term is a $dx \cup dt_\xi$
term for each tile $\xi$ for which $\rho(\xi)$ is a word beginning with
$\zeta$, plus a $dy \cup dt_\xi$ term for each tile $\xi$ for which
$\rho(\xi)$ contains $\zeta$ in the middle or end of the word.
If the substitution forces the border in $m$ steps, then $\rho^m$ of  
each
edge emerging from a vertex in $P$ is a word beginning with the same
letter (call it $\omega$).  $(\rho^m)^* (\mu_x \cup dt_\omega)$ then
equals the sum of the $\mu_x \cup dt$ terms from all the edges emerging
from this vertex (plus additional $\mu_y \cup dt$ terms).  However, the
sum of the $\mu_x \cup dt$ terms from the edges emerging from
the vertex equals the sum of the $\mu_y \cup dt$ terms from the edges
entering the vertex.  Likewise, if $\eta \ne \omega$, then $(\rho^m)^*
(\mu_x \cup dt_\eta)$ contains none of the $\mu_x \cup dt$ terms from
edges emerging from the vertex.  Either way, the pullback by $\rho^m$ of
a $\mu_x \cup dt$ term can be expressed as a sum of $\mu_y \cup dt$  
terms.
The same argument (applied to ends of words) shows that the pullback
by $\rho^m$ of a $\mu_z \cup dt$ term can be expressed as a sum of
$\mu_y \cup dt$ terms.  In particular, the $\mu_x \cup dt$ and $\mu_z  
\cup
dt$ terms that are linearly independent of the $\mu_y \cup dt$ terms  
do not
appear in the eventual range of $\rho^*$, and hence do not contribute to
the direct limit of $H^2(L)$.

In dimension 3, $\rho^*$ sends $\mu_x\cup\mu_z \cup dt$ terms to
$\mu_x \cup \mu_y \cup dt$ and $\mu_y \cup \mu_z \cup dt$ terms.  
Therefore,
we need only consider the direct limit of $\mu_x \cup \mu_y \cup dt_ 
\zeta$
and $\mu_y \cup \mu_z \cup dt_\zeta$ terms. We associate
$\mu_x \cup \mu_y \cup dt_\zeta$ with the vertex in $P$ that $\zeta$  
leads
into, and associate $\mu_y \cup \mu_z \cup dt_\zeta$ with the vertex  
it leads
out of.  As in dimension 2, forcing the border implies that
$(\rho^m)^*$ takes each $\mu_x \cup \mu_y \cup dt$
term to either all or none of the $\mu_x \cup \mu_y \cup dt$ terms  
associated
with a vertex.  When dealing with $\mu_x \cup \mu_y \cup dt$, we may  
therefore
restrict our attention to sums of all the edges emerging from a vertex,
and when
dealing with $\mu_y \cup \mu_z \cup dt$ we may restrict our attention  
to sums
of all the edges leading into a vertex.  However, the sum of all the $ 
\mu_y
\cup \mu_z \cup dt$s from the edges leading into a vertex equals the sum
of all the $\mu_x \cup \mu_y \cup dt$s leading out of the vertex, so we
have exactly one independent term per vertex.

The substitution $\rho$ maps vertices of $P$ to vertices, and has a
natural pullback action on the
3-forms associated to vertices.  Eventually, $\rho$ merely permutes the
vertices that describe boundaries of two infinite-order supertiles.  The
3-forms associated with these vertices give a basis for $H^3(X_\Sigma)$.

In dimension 4, the pullback map is zero, so $H^4(X_\Sigma)=0$.

Finally, we consider the coefficient of $\mu$ in $H^2(X_\Sigma)$ and the
coefficient of $\mu\otimes\mu$ in $H^3(X_\Sigma)$.  The first is the  
direct
limit of the transpose of the substitution matrix as applied to edges  
of $P$,
and the second is the direct limit of the transpose of the substitution
matrix as applied to vertices of $P$.  In other words, they are the
direct limit of 1-cochains and 0-cochains on $P$ under substitution.
These are closely
related to direct limits of the cohomology of $P$ (i.e., to the
cohomology of $X_\rho$).  Since the direct limit of the
0-cochains (namely $\Z^n$) is $\Z^{n-1}$ more
than $H^0(X_\rho)$, the direct limit of the 1-cochains
must be $\Z^{n-1}$ more than $H^1(X_\rho)$. \qed

In stating Theorem 1, we assumed that every boundary between infinite- 
order
supertiles in $X_\rho$ led to a fault line.  As the following example  
shows,
this is not always the case.

Let $A_1$, $A_2$ and $A_3$
be tiles of width $\lambda$ and height 1/3, and let $B_1$, $B_2$, and  
$B_3$
be tiles of width 3 and height 1/3.
Our substitution is
\begin{eqnarray}
\Sigma(A_1) = \begin{pmatrix} B_2 & A_2 \cr B_1 & A_1 \end{pmatrix} &
\Sigma(A_2) = \begin{pmatrix} A_1 & B_1 \cr B_3 & A_3 \end{pmatrix} &
\Sigma(A_3) = \begin{pmatrix} A_3 & B_3 \cr A_2 & B_2 \end{pmatrix} \\
\notag \Sigma(B_1) = \begin{pmatrix} A_2 & A_2 & A_2 \cr A_1 & A_1 & A_1
\end{pmatrix} &
\Sigma(B_2) = \begin{pmatrix} A_1 & A_1 & A_1 \cr A_3 & A_3 & A_3
\end{pmatrix} &
\Sigma(B_3) = \begin{pmatrix} A_3 & A_3 & A_3 \cr A_2 & A_2 & A_2
\end{pmatrix}.
\end{eqnarray}
Each row in a tiling contains either $A_1$s and $B_1$s (call
this an $\alpha$ row), $A_2$s and $B_2$s ($\beta$) or $A_3$s and $B_3$s.
The complex $P$ consists of three edges ($\alpha$, $\beta$,
$\gamma$) running in a circle, with $\alpha$ followed by $\beta$  
followed
by $\gamma$ followed by $\alpha$.  The vertical substition $\rho$ is
\begin{equation} \rho(\alpha) = \alpha \beta, \qquad
\rho(\beta) = \gamma \alpha, \qquad \rho(\gamma)=\beta \gamma.
\end{equation}
Fault lines do not develop at the
boundary of $\alpha$ and $\beta$ supertiles, or at the boundary of
$\beta$ and $\gamma$ supertiles, just at the boundary of $\gamma$ and
$\alpha$.

\begin{figure}
\vbox{\epsfxsize=2truein\epsfbox{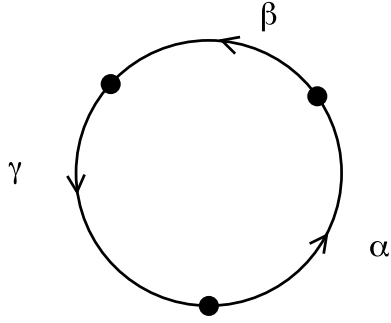}}
\caption{The Anderson-Putnam complex for the vertical tiling space}
\end{figure}

This example may look mysterious, but it is just a rewriting of
the basic example of section \ref{simple}.
$A_1$, $A_2$, and $A_3$ are just the bottom,
middle and top thirds of the $A$ tile, while $B_1$, $B_2$ and $B_3$
are the bottom, middle and top thirds of the $B$ tile.  The  
cohomology is
correctly described by Theorem 1, only with $n$ being the number of  
possible
fault lines (namely 1), not the number of ways that two infinite-order
supertiles can meet (namely 2).  This observation generalizes to

\begin{thm}
\label{thm2}
   Suppose we have a 2-dimensional substitution $\Sigma$,
   generated by a vertical substitution $\rho$ and
   horizontal substitutions $\sigma_1, \ldots, \sigma_N$, as in  
Theorem 1.
   Suppose that the boundaries
   between infinite-order $\rho$-supertiles are either fault lines or
   are rigid, with the patterns on both sides of the boundary being
   mutually locally derivable.  As before, let $\mu=H^1(X_\sigma)$ and
   let $\nu=H^1(X_\rho)$.  Let $n$ be the number of ways that a fault
   line can develop in an $X_\Sigma$ tiling.  Then the cohomology of
   $X_\Sigma$ is as follows:
\begin{eqnarray} \notag H^0(X_\Sigma)& = & \Z, \\
\notag H^1(X_\Sigma) & = & \nu, \\
H^2(X_\Sigma) & = & \mu \otimes (\nu \oplus \Z^{n-1}), \\
\notag H^3(X_\Sigma) & = & (\mu \otimes \mu)^n, \\
\notag H^k(X_\Sigma) & = & 0 \hbox{ for }k>3.
\end{eqnarray}
\end{thm}

\noindent Proof: Let $P$ be the Anderson-Putnam complex of $\rho$.
By assumption, each vertex of $P$ either generates a fault line
(call this an essential
vertex) or has the patterns on both sides of the vertex precisely  
aligned.
We rewrite $\rho$ using the essential vertices as stopping and  
starting rules as in \cite{BDstartstop},
and rewrite $\Sigma$ in terms of these new vertical tiles.  By  
construction,
each of the $n$ vertices of the new vertical substitution generates a  
fault
line, so Theorem 1 applies directly. \qed

\section{Open problems}
\label{open}

\begin{enumerate}
\begin{item}
We understand fault lines for substitutions on two letters, but what
about more complicated substitutions?  Is it possible to have lines
without finite local complexity that do not allow arbitrary shears?
What can we say about the cohomology of tiling spaces that allow such
``irregular'' fault lines?
\end{item}\begin{item}
In considering Theorem \ref{thm2}, we assumed that
all boundaries between infinite-order $\rho$-supertiles either generated
fault lines or had the two sides remain in lockstep.  When the  
horizontal
stretching factor is not Pisot, are these
the only possibilities?
\end{item}\begin{item}
What happens if the horizontal substitutions of Theorem 1 are different
enough that the rows in a supertile are not all the same
up to translation?  (The horizontal stretching factors would all have to
be the same, which would constrain the possible abelianizations of the
different $\sigma_i$, but the actual substitutions could differ.)
There are natural conjectures for what $H^1$ and
$H^3$ of such
a tiling space should look like.  $H^1$ should come entirely from the
vertical substitution and $H^3$ should contain a copy of $\mu_1
\otimes \mu_2$ for each possible fault line, where $\mu_1$ (resp. $ 
\mu_2$)
is $H^1$ of the tiling space that describes the row immediately above
(below) the fault line.   However, it is not at all clear what $H^2$  
should
be.
\end{item}\begin{item}
Some tilings (for instance the one in Figure 2) allow both vertical and
horizontal fault lines.
Since a single tiling can exhibit a lack of coordination across
at most one fault line, it is
easy to guess what $H^3$ of such a tiling space should look like, with a
contribution (as in the previous problem)
from each possible fault line, vertical or horizontal.  But does $H^1$
vanish altogether?  What about $H^2$?
\end{item}\begin{item}
Up to this point we have only been considering rectangular tiles. What
if the tiles are not rectangular, as in the generalized pinwheel? Again,
it is possible to predict the cohomology in the highest dimension (4  
for the
generalized pinwheel, since the rotations provide an additional
degree of freedom), with a contribution from each possible species of
fault line.  However, we do not venture to guess
the lower dimensional cohomology.
\end{item}
\end{enumerate}

\section{Acknowledgements}

Much of this work was done at the Banff International Research Station
in 2005.  All of the participants in the Focused Research Group on
Topological Methods in Aperiodic Tilings contributed substantially to  
the ideas
in this paper.  In particular we would like to thank Marcy Barge,  
Beverly Diamond,
John Hunton, Johannes Kellendonk, and Ian Putnam.  Additional work  
was done at the Aspen Center
for Physics in 2006.  The work of Lorenzo Sadun is partially  
supported by
the National Science Foundation under grant DMS-0401655.

\end{document}